\pdfoutput=1
\documentclass[a4paper,11pt]{article}
\usepackage{amsfonts,latexsym,rawfonts,amsmath,amssymb,amsthm}
\usepackage{cases, cite}
\usepackage{hyperref}

\hypersetup{
    hidelinks,
}
\usepackage{graphicx}
\usepackage{setspace}
\usepackage{titlesec}
\usepackage{enumitem}
\usepackage{soul}

\usepackage{tikz}
\usepackage{scalerel}
\usepackage{pict2e}
\usepackage{tkz-euclide}
\usetikzlibrary{calc}
\usetikzlibrary{patterns,arrows.meta}
\usetikzlibrary{shadows}
\usetikzlibrary{external}
\usetikzlibrary{patterns.meta}

\usepackage{caption}
\usepackage{subcaption}

\usepackage{nicefrac}

\titleformat{\chapter}[block]
    {\normalfont\huge\bfseries}{\thechapter}{20pt}{\Huge}

\RequirePackage{color}

\theoremstyle{plain}
  \newtheorem{theorem}{Theorem}[section]
  \newtheorem{lemma}[theorem]{Lemma}
  \newtheorem{proposition}[theorem]{Proposition}
  
  \newtheorem{corollary}[theorem]{Corollary}

  \newtheorem{remark}[theorem]{Remark}
  \newtheorem{definition}[theorem]{Definition}

\theoremstyle{definition}

\numberwithin{equation}{section}

\DeclareMathOperator*{\vol}{vol}
\DeclareMathOperator*{\intr}{int}

\newcommand{\R}{\mathbb R}

\def\sub{\subseteq}%

\def\RR{\mathbb{R}}%

\usepackage[nobysame,initials]{amsrefs}
\usepackage{xcolor}

\let\originalleft\left
\let\originalright\right
\renewcommand{\left}{\mathopen{}\mathclose\bgroup\originalleft}
\renewcommand{\right}{\aftergroup\egroup\originalright}

\usepackage[utf8]{inputenc}
\title{Boundary restricted Brunn-Minkowski inequalities \thanks{The authors were partially supported by ISF grant number 784/20. The first and second named authors were also partially supported by European Research Council (ERC) under the European Union's Horizon 2020 research and innovation programme (grant agreement No. 770127).}}

\author{Shiri Artstein-Avidan\thanks{School of Mathematics, Tel Aviv University, Tel Aviv, Israel; \href{mailto:shiri@tauex.tau.ac.il}{shiri@tauex.tau.ac.il}.}  \and Tomer Falah\thanks{School of Mathematics, Tel Aviv University, Tel Aviv, Israel; \href{mailto:tomerfalah@mail.tau.ac.il}{tomerfalah@mail.tau.ac.il}.} \and Boaz A. Slomka\thanks{Department of Mathematics, The Open University of Israel, Ra'anana, Israel; \href{mailto:slomka@openu.ac.il}{slomka@openu.ac.il}.}}
\date{}

\begin{document}

  \maketitle

\begin{abstract}
In this paper  we explore questions regarding the Minkowski sum of  the boundaries of convex sets. 
Motivated by a question suggested to us by V.~Milman regarding the volume of $\partial K+ \partial T$ where $K$ and $T$ are convex bodies, we prove sharp volumetric lower bounds for the Minkowski average of the  boundaries of sets with connected boundary, as well as some related results. 
\end{abstract}



\section{Introduction}
\renewcommand{\thefootnote}{\fnsymbol{footnote}} 
\footnotetext{\emph{Key words and phrases.} Brunn-Minkowski inequality, convex bodies, polytopes, zonotopes, covering.}
\renewcommand{\thefootnote}{\arabic{footnote}} 
\renewcommand{\thefootnote}{\fnsymbol{footnote}} 
\footnotetext{\emph{2010 Mathematics Subject Classification.} 52A40, 52B12, 52C17.}
\renewcommand{\thefootnote}{\arabic{footnote}} 
%
The classical Brunn-Minkowski inequality (see e.g. \cite{Sch_book}*{Section 6.1}, \cite{AGM}*{Section 1.2}) states that 
\begin{equation}\label{eq:BMadditive} \vol(A+B)^{1/n}\ge\vol(A)^{1/n}+\vol(A)^{1/n}\end{equation}
for any $A$ and $B$ which are non-empty Borell subsets of $\R^n$. Here $\vol(\cdot)$ denotes the usual Lebesgue measure on $\RR^n$ and $A+ B = \{ a+b : a\in A, b\in B\}$. 
When $A,B$ are convex bodies (compact convex sets with non-empty interior), equality holds in \eqref{eq:BMadditive} if and only if some translate of $A$ is homothetic to $B$.
An equivalent geometric form of this inequality is the statement that for every  $\lambda \in [0,1]$
\begin{equation}\label{eq:BMgeom}\vol(\lambda A+(1-\lambda )B)\ge\vol(A)^\lambda\vol(B)^{1-\lambda}.\end{equation}
When $A,B$ are convex bodies and $\lambda \in (0,1)$ equality holds in \eqref{eq:BMgeom} if and only if $A$ is a translate of $B$.

Inspired by results on Brunn-Minkowski type inequalities for the lattice point enumerator (see e.g.~\cites{HCLY22 ,KL19, IYNZ, ILYN22, AGLYN,HKS}) where an error term emerges due  to lattice points near the boundary of a set, Prof.~V.~Milman asked whether a complementary phenomenon might hold, namely whether discarding all interior points of a convex body, and taking the Minkowski sum of two boundaries, one might retrieve a substantial amount of the sum of the bodies.

In this note we show that indeed  a Brunn-Minkowski inequality {\em of geometric form} holds in this case for the boundary average of two convex bodies in $\RR^n$ starting from dimension $2$ (in dimension $n=1$ the boundary sum is simply $3$ or $4$ points). 
In fact, as is the case in the classical Brunn-Minkowski inequality, one can relax the convexity condition. We will assume only that each of the sets is compact and with a connected boundary. Our main theorem is as follows.

\begin{theorem}\label{thm:av}
Let $n\ge 2$ and let $K,T\subseteq \RR^n$  be  compact sets with connected boundary. It holds that 
	\[ \vol\left(\frac{\partial K + \partial T }{2}\right) \ge \sqrt{\vol(K) \vol(T)}.\]
When $K$ and $T$ are convex bodies, equality holds if and only if either $T$ is a translation of $K$, or $n=2$ and  $K,T\sub\RR^2$ are homothetic and centrally-symmetric up to translation.
\end{theorem}

Theorem \ref{thm:av} directly extends to multiple boundary sums, as follows.
\begin{corollary}\label{cor:multiBBM}
Let $n\ge 2$ and let $K_{1},\dots,K_{m}\sub\R^{n}$ be  compact sets with connected boundary, where $m>2$. It holds that
\[
\vol\left(\frac{\partial K_{1}+\dots+\partial K_{m}}{m}\right)\ge\vol(K_{1})^{\frac{1}{m}}\dots\vol(K_{m})^{\frac{1}{m}}.
\]
When $K_1,\dots,K_m$ are convex bodies, equality holds if and only if  $K_1,\dots,K_m$ are  translations of one another. 
\end{corollary}

Replacing $\lambda = 1/2$ in Theorem \ref{thm:av} with some other $\lambda\in (0,1)$, one cannot expect a similar analogue of \eqref{eq:BMgeom}, simply by considering two copies of the unit ball (similarly, a result of the form \eqref{eq:BMadditive} cannot hold for boundary sums by considering two balls of different radii). Note that this does not contradict 
Corollary \ref{cor:multiBBM}, say for $\lambda =\frac13$, since in general $\partial K+\partial K +\partial T \neq 2\partial K+\partial T$.

Nevertheless, we show the following theorem for general $\lambda\in (0,1)$, the case $\lambda=\frac 12$ of which is Theorem \ref{thm:av}.
\begin{theorem}\label{thm:BBM}
Let $n\ge 2$, and let  $K,T\subseteq \RR^n$ be compact sets with connected boundary. For any $\lambda \in(0,1)$ it holds that 
	$$
	\vol(\lambda\partial K+(1-\lambda)\partial T)\vol(\lambda\partial T+(1-\lambda)\partial K)\ge \vol(K)\vol(T)\cdot (1-|1-2\lambda|^n)^2.
	$$
When $K$ and $T$ are convex bodies and $\lambda\neq1/2$, equality holds if and only if either $K$ and $T$ are  translations of some centrally-symmetric body, or $n=2$ and $K,T\sub\RR^2$ are translations of homothets of some centrally-symmetric body.


\end{theorem}

In the final part of the paper we prove another Brunn-Minkowski-type inequality, this time of arithmetic type, for boundary sums. To motivate this, let us note that our inequalities can be viewed within the framework of ``restricted Minkowski sums'' a la Szarek-Voiculescou. More precisely, given two sets $A,B\sub\RR^n$ and a subset $\Theta \sub A\times B$, one defines their $\Theta$-restricted Minkowski sum
by 
\[ A +_\Theta B = \{ x+y: (x,y)\in \Theta\}.
\]

In our case, of course,  $\Theta = \partial A \times \partial B$ which is ``geometrically large'' in some sense, but has zero volume. Szarek and Voiculescu \cite{SV} proved  a complementing result which is  that if $\vol(\Theta)$ is sufficiently close to $\vol(A\times B)$ then  an (arithmetic)  Brunn-Minkowski inequality with exponent $2/n$ holds. We quote their theorem as Theorem \ref{thm:restriceted sums} below.
It turns out that utilizing their theorem we may provide a corresponding version of Brunn-Minkowski inequality for our boundary sums, under a volume ratio restriction. This is carried out in Theorem \ref{thm:restricted_BM}.

The note is organized as follows: In  Section \ref{sec:deomp} we discuss the geometric structure of boundary sums. In Section \ref{sec:volbd} we prove our main geometric Brunn-Minkowski type theorems for boundary sums. Finally, in Section \ref{sec:restrict} we discuss restricted Minkowski sums, state and prove our arithmetic version of Brunn-Minkowski inequality for boundary sums. 



\section{Decomposition of the sum}\label{sec:deomp}

We shall need the following definition of the ``Minkowski difference'' of two sets  $K \ominus T$, which is similar, but not identical, to the definition of ``Minkowski subtraction'' appearing in \cite{Sch_book}*{Section 3.1}. 
\begin{definition}\label{def:MinkDif}
	Let $K,T\subseteq \RR^n$. Their \emph{Minkowski difference} is given by the set 
	\[
	K \ominus T = \bigcap_{t \in T} (t + \intr(K))
 = \{x: x-T \subseteq \intr(K)\}.
	\]
\end{definition}

\noindent  In other words, $K\ominus T$ is the largest (with respect to inclusion)  set $L$ such that $L-T \sub {\rm int}(K)$. 
For example, for the cube $Q^n = [-2,2]^n$ of side length $4$ and the Euclidean unit ball $B_2^n$  (see Figure \ref{fig:Set diffrnce})
\[ Q^n \ominus B_2^n = 
\{ x: x+B_2^n \sub (-2, 2)^n\} = (-1,1)^n = Q^n \ominus [-1,1]^n.
\] 







\begin{figure}[h]
    \centering





\begin{tikzpicture}[scale = 0.7][!h]

\draw (3,0) circle[radius=1] node at (3,-1.5) {$B_2^2$};

\draw[pattern={Lines[angle=45]}] (5,-2) rectangle (9,2) node at (7,-2.5) {$Q^2$};

\draw  (11,-2) rectangle (15,2) ;
\draw[dashed] (11,-1) arc (180:270:1) -- (14,-2) arc (-90:0:1)--(15,1) arc (0:90:1)--(12,2) arc (90:180:1) -- cycle;
\node at (13,-1.5) {$Q^2\ominus B_2^2$};
\draw[dashed, pattern={Lines[angle=45]}] (12,-1) rectangle (14,1);
\end{tikzpicture}

    \caption{Minkowski {difference} of a square and a disk, both centered at the origin.}
    \label{fig:Set diffrnce}
\end{figure}
\begin{remark}
Before moving on, let us mention a few curious facts. Clearly for a centrally symmetric convex body $K$ and $t\in (0,1)$  we have that $K \ominus (tK) = (1-t){\rm int}(K)$. It is not hard to check that this is no longer true for non-symmetric $K$. Sightly less obvious is the fact that for a simplex $\Delta$ we have that $\Delta \ominus K$ is always either empty or a dilation of ${\rm int }\Delta$ (up to translation). For example, $\Delta \ominus \frac13 \Delta = \frac13 {\rm int} \Delta$, see Figure \ref{fig:Simplex_diff}. In dimension $n=2$, the closure of the difference body $K\ominus T$ is always a {\em summand} of 
$K$, see \cite{Sch_book}*{Theorem 3.2.4}.\end{remark}







\begin{figure}[ht]
    \centering

\begin{tikzpicture}[scale = 1][!h]

\draw 
     (0,0) 
  -- (1,0) 
  -- (0,1) node at (0.5, -0.5){$\frac13\Delta$}
  -- cycle;

\draw[pattern={Lines[angle=45]}] 
     (2,0)   
  -- (5,0) 
  -- (2,3) node at (3.5,-0.5){$\Delta$}
  -- cycle;

\draw 
     (9,0) 
  -- (6,0) 
  -- (6,3) ; 
\draw 
    (6,3) 
   -- (7,2);

\draw
      (8,1) 
   -- (9,0);

 \draw [dashed,pattern={Lines[angle=45]}]
     (7,1) 
  -- (8,1) 
  -- (7,2) node at (7.5,0.5)
  {$\Delta\ominus\frac13\Delta$}
  -- cycle;


\end{tikzpicture}

    \caption{Minkowski {difference} of homothetic simplices $\Delta$ and $\frac13\Delta$} (where the lower-left vertex of each simplex is at the origin).
    \label{fig:Simplex_diff}
\end{figure}

\noindent Our first observation is that the Minkowski sum of $K$ and $T$ decomposes   as follows.
\begin{proposition}\label{prop:decomp-of-sum}
Let $n\ge2$, and let $K,T\sub\RR^n$ be compact sets, each with connected boundary. Assume $\vol(K) \ge \vol(T)$. It holds that  $(T\ominus K)= \emptyset$ and
\begin{equation}\label{eq:decomposition}
	K + T  = K + \partial T = (\partial K + \partial T ) \cup (K\ominus T). \end{equation}
 Moreover, the union is disjoint, and it is also true that 
 \begin{equation}\label{eq:decomposition2}
	  \partial K + \partial T    = \partial K + T . \end{equation}
\end{proposition}
\noindent For an illustration of this decomposition in the case of a square and a disk, see Figure~\ref{fig:Partial sum}. 






\begin{figure}[h]
    \centering

\begin{tikzpicture}[scale = 0.7][!h]

\draw[pattern={Lines[angle=45]}] (3,0) circle[radius=1] node at (3,0-1.5) {$T$};

\draw (5,-2) rectangle (9,2) node at (7,-2.5) {$K$};

\draw[pattern={Lines[angle=45]}] (10,-2) arc (180:270:1) -- (15,-3) arc (-90:0:1)--(16,2) arc (0:90:1)--(11,3) arc (90:180:1) -- cycle;
\node at (13,-3.5) {$\partial K+\partial T$};
\fill[white] (12,-1) rectangle (14,1);
\draw[thick] (12,-1) rectangle (14,1);
\end{tikzpicture}

    \caption{Boundary sum of a square and a disk.}
    \label{fig:Partial sum}
\end{figure}

\begin{remark} \label{minus}
Note that when no translate of $-K$ can be shifted into ${\rm int}(T)$, and no translate of $-T$ can be shifted into ${\rm int}(K)$, we have  $\partial K + \partial T  = K+T$. This fact was noticed already for example in \cite{Ulivelli23}. 
	This is the case, for example, if in some direction the width of $K$ exceeds than the width of $T$, and in another direction the width of $T$  exceeds the width of $K$.  In the case of $T=K$, the fact that $\partial K+ \partial K = K + K$ was noted and put to use in \cite{FLZ22}. Iterating \eqref{eq:decomposition}, we get that the Minkowski sum of $k$ copies of $\partial K$  equals to the sum of $k$ copies of $K$. 
\end{remark}

To prove Proposition \ref{prop:decomp-of-sum} we shall use the following intuitive lemma. 

\begin{lemma}\label{lem:bc}
   Let $n\ge2$, and let $K,T\sub\R^{n}$ be compact sets, each with connected boundary. If $\partial T\subseteq \intr(K)$ then $T\subseteq \intr(K)$. 
\end{lemma}
\begin{proof}
Assume $\partial T\subseteq \intr(K)$. 
The set $\R^n\setminus \partial T$ is open, so we may write it as a disjoint union of
its connected components which are themselves open sets. Denote these sets by $\{C_{i}\}_{i\in I}$. 
Let $J\sub I$ consist of those $C_i$ which intersect $\intr (T)$, and then $\intr(T) = \cup_{j\in J} C_j$. 
Since $\{C_{i}\}_{i\in I}$ are disjoint and open it follows that 
\[\partial C_i \subseteq \R^n\setminus\cup_{m\in I}C_{m}= \partial T.\]

The fact that $\partial T\subseteq \intr(K)$ implies that $\partial T\cap\partial K=\emptyset$ and thus  
$
\partial K\subseteq\cup_{i\in I}C_{i}$. 
Since $\partial K$ is connected there exists a $k\in I$ such that $\partial K\subseteq C_{k}$. 
Let us show that for $i\ne k$,  $\intr(K)\supseteq C_{i}$. By disjointness, clearly, 
either $\intr(K)\supseteq C_{i}$ or $\intr(K)\cap C_{i}=\emptyset$.
However, since $\partial C_i \subseteq \intr(K)$ it must be the case that $C_{i}\cap \intr(K)\ne\emptyset$, implying $C_i\subseteq \intr(K)$
for any $i\ne k$.
As $T$ and $K$ compact we know that $\R^n\setminus T$ is unbounded, and  $\intr(K)$ is bounded, meaning $\R^n\setminus T \nsubseteq \intr(K)$. 
In particular, $\cup_{j\not\in J}C_j =  \R^n\setminus T\nsubseteq \cup_{i\neq k} C_i$. Thus $k\not\in J$, that is, $C_k \subseteq \R^n\setminus T$.
We get that 
\[T \subseteq \partial T \cup \bigcup_{i\in J}C_{i} \subseteq \intr(K).\]
\end{proof}
\begin{proof}[Proof of Proposition \ref{prop:decomp-of-sum}]
The fact that $T\ominus K = \emptyset$ is immediate due to the volume assumption. Indeed,  $x-K$ is not included in  $\intr (T)$ for any $x$, as the complement $\intr (T) \setminus (x-K)$ would then be a non-empty open set with zero volume. 

We start by proving the  decomposition 
\begin{equation}\label{eq:decomp111} (\partial K +  T ) \cup (K\ominus T) = K + T .\end{equation}

Note that this union is disjoint. Indeed, by definition,  for any $x\in K\ominus T$, we have $x-T\subseteq \intr(K)$ and,  in particular,  $(x-T)\cap \partial K=\emptyset$, implying that  the sets are disjoint.

To prove \eqref{eq:decomp111}, we show a double inclusion. Clearly $\partial K +  T  \subseteq K + T $, and 
 $K\ominus T \subseteq  K + T $ since 
if $x-T\sub\intr (K)$ then $x\in T+K$. 

For the opposite inclusion,  let $x \in (K + T)\setminus(K \ominus T)$ (we should show $x\in \partial K + T$). 
Let  $t_1 \in T$ such that $x-t_1 \in K$. Since $x \notin K \ominus T$, there exist $t_2 \in T$ such that $x-t_2\not\in  \intr(K)$.
    If either $x-t_1\in\partial K$ or $x-t_2\in\partial K$ then $x\in \partial K + T$ as needed. Otherwise, $x-T$ intersects both $\intr(K)$ and $\R^n\setminus K$ which are open disjoint sets. As $x-T$ is   connected, it follows that it must intersect their complement $\partial K$,  
 as claimed. This proves \eqref{eq:decomp111}.

 To prove \eqref{eq:decomposition}, we note that   the inclusions
\[ 
	K + T  \supseteq  K + \partial T \supseteq (\partial K + \partial T ) \cup (K\ominus T) \]
  are immediate. 
  For the less trivial  inclusion 
 $K + T  \subseteq (\partial K + \partial T ) \cup  (K\ominus T)$
 we use \eqref{eq:decomp111}, so it is enough to show 
    \[
    \partial K +  T \subseteq \partial K + \partial T .
    \]
  (Note that this inclusion is in fact an 
   equality here since the opposite inclusion is trivial). 
  
  Let $x \in \partial K + T$. Let  $k_1 \in \partial K$ such that $x-k_1 \in T$.
Since $T\ominus K = \emptyset$ we know that 
 $x-K\not\sub\intr(T)$. 
By Lemma \ref{lem:bc} this means that  $x-\partial K\not\sub \intr(T)$.
Therefore, there is some    
   $k_2 \in \partial K$ such that $x-k_2 \not\in \intr(T)$.
 If either $x-k_1$ or $x-k_2$ are in $\partial T$ then $x\in \partial K + \partial T$ as needed. 
  Otherwise $x-k_1 \in \intr(T)$ and $x-k_2 \in \R^n \setminus T$.
   Since $\intr(T)$ and $\R^n\setminus T$ are open disjoint sets   intersecting the connected set $(x-\partial K)$,   it must be the case that $(x-\partial K)$  intersects $\partial T$, as needed. 
\end{proof}

\section{Volume bounds}\label{sec:volbd}

Although Theorem \ref{thm:av} is a special case of Theorem \ref{thm:BBM}, it admits a ``one line'' proof, as we next present. 
\begin{proof}[Proof of Theorem \ref{thm:av}]
From Proposition \ref{prop:decomp-of-sum} we know that if $K\ominus T = T\ominus K = \emptyset$ then $\partial K + \partial T = K+T$ and the claim follows from the classical Brunn-Minkowski inequality \eqref{eq:BMgeom}. 
We may thus assume without loss of generality that $K\ominus T \ne \emptyset$, and by Remark \ref{minus}  we know then that $T\ominus K = \emptyset$. Denote $L =K\ominus T$. 
Since $L-T \sub\intr(K)$, the Brunn-Minkowski inequality \eqref{eq:BMadditive} implies $\vol(K)^{1/n} \ge \vol(L)^{1/n} + \vol(T)^{1/n}$.  Using this,  Proposition \ref{prop:decomp-of-sum}, \eqref{eq:BMadditive} and that $(c+d)^k \ge c^k + d^k$ for positive $c,d$ and $k\ge 1$ we get: 
	\begin{eqnarray*}
		\vol\left(\frac{\partial K + \partial T}{2}\right) & = & \frac{1}{2^n}\left(\vol(K + T) - \vol(L)\right)\\
		& \ge & \left( \frac{\vol(K)^{1/n} + \vol(T)^{1/n}}{2}\right)^n - \frac{ \vol(L)}{2^n} \\
		&= & \left(\left( \frac{\vol(K)^{1/n}+ \vol(T)^{1/n}}{2}\right)^2\right)^{n/2} - \frac{ \vol(L)}{2^n}\\
 	&= & \left(\left( \frac{\vol(K)^{1/n}-\vol(T)^{1/n}}{2}\right)^2
 	+\vol(K)^{1/n}\vol(T)^{1/n}
 	\right)^{n/2} - \frac{    \vol(L)   }{2^n}\\
  & \ge & \left( \frac14\vol(L)^{2/n}
  +\vol(K)^{1/n}\vol(T)^{1/n}
  \right)^{n/2} - \frac{ \vol(L)}{2^n} \\
  & \ge &  \frac{\vol(L)}{2^n} 
  +\vol(K)^{1/2}\vol(T)^{1/2}
    - \frac{ \vol(L)}{2^n}  = \vol(K)^{1/2}\vol(T)^{1/2},
	\end{eqnarray*}
	as required. Note that we assumed $n\ge 2$ (so that the power $n/2\ge 1$).

   To analyze the equality case note that we used the classical Brunn-Minkowski twice. First to say $\vol(K+T)\ge (\vol(K)^{1/n}+\vol(T)^{1/n})^n$,  equality for which is when $K$ and $T$ are homothetic (up tp translation). The second  use of  Brunn-Minkowski is for $\vol(K)^{1/n}-\vol(T)^{1/n}\ge\vol(L)^{1/n}$, for which the equality cases are as follows. 
   One option is that $\vol(K)=\vol(T)$ (and $L$ is then empty) which in turn means that $T$ is a translation of $K$. 
   The other option is that $L-T=K$ and (as when applying Brunn-Minkowski to this sum we get equality) $L$, $K$ and $-T$ are homothetic to one another (and to $T$, up to taking the closure of $L$, as it is defined not to be closed). This means in turn that $T$ is, up to translation, centrally symmetric (as well as its homothet $K$ of course). In addition,  for the inequality $(c+d)^k \ge c^k + d^k$ to be satisfied as an equality, either  $k=1$ (and $c,d$ are any numbers) or, if $k >1$, one must have that $c=0$ or $d=0$, which is the case if any only if $L=\emptyset$, namely (as we already showed that $K$ and $T$ must be homothetic) $T$ is a translation of $K$, meaning the only equality case where $K,T$ are not translations of each other is for $n=2$ and $K,T$ being homothehic to each other and centrally-symmetric up to translation. \end{proof}

To prove Corollary \ref{cor:multiBBM}, we need the following simple lemma.
\begin{lemma}
\label{lem:pBM_balls}Let $m\ge1$. For any $x_{1},\dots,x_{m+1}\ge0$
such that $x_{1}+\dots+x_{m}\le x_{m+1}$,
\[
\left(\frac{2}{m+1}\right)\sqrt{\left(x_{1}+\dots+x_{m}\right)\, x_{m+1}}\ge x_{1}^{\frac{1}
{{m+1}}}\dots x_{m+1}^{\frac{1}{m+1}}.
\]
For $x_{m+1}>0$ and $m>1$, the inequality is strict. 
\end{lemma}

\begin{proof}
By homogeneity of the inequality, it is enough to show that for all non-negative $x_{1},\dots,x_{m}$ such
that $x_{1}+\dots+x_{m}\le1$, one has
\[
\left(\frac{2}{m+1}\right)\left(x_{1}+\dots+x_{m}\right)^{\frac{1}{2}}\ge x_{1}^{\frac{1}{{m+1}}}\dots x_{m}^{\frac{1}{m+1}}.
\]
By the arithmetic-geometric means inequality, it suffices to show
that 
\[
\left(\frac{2}{m+1}\right)^{2\left(m+1\right)}\left(x_{1}+\dots+x_{m}\right)^{m+1}\ge\left(\frac{x_{1}+\dots+x_{m}}{m}\right)^{2m}
\]
or, equivalently, that 
\[
\frac{1}{m}\left(\frac{2m}{m+1}\right)^{\left(m+1\right)}\ge\left(x_{1}+\dots+x_{m}\right)^{\frac{m-1}{2}}.
\]
Since $1\ge x_{1}+\dots+x_{m}$ and $f\left(m\right)=\frac{1}{m}\left(\frac{2m}{m+1}\right)^{m+1}\ge f\left(1\right)=1$
is easy checked to be a strictly monotonically increasing function of $m\ge 1$, the assertion of the lemma follows.
\end{proof}

\begin{proof}[Proof of Corollary \ref{cor:multiBBM}]
By \eqref{eq:decomposition2}, either $\partial K_{i}+\partial K_{j}=K_{i}+\partial K_{j}$
or $\partial K_{i}+\partial K_{j}=\partial K_{i}+K_{j}$ for each
$i\neq j$. Hence, by relabeling the sets if needed, we may assume without loss of generality
that $\partial K_{1}+\dots+\partial K_{m}=K_{1}+\dots+K_{m-1}+\partial K_{m}$.

If for any $x\in\RR^n$, $x+K_{1}+\dots+K_{m-1}\not\sub K_{m}$   then $K_{1}+\dots+K_{m-1}+\partial K_m=K_{1}+\dots+K_{m}$
and the proposition follows from the classical Brunn-Minkowski inequality (with the equality cases only when the bodies are translates of one another).

Next, we assume we are in the complementary case, namely that $K_{1}+\dots+K_{m-1}\sub K_{m}$ (assuming $x = 0$ without loss of generality). Denote $x_{i}=\vol(K_{i})^{1/n}.$
Note that, by Brunn-Minkowski inequality, $x_{1}+\dots +x_{m-1}\le x_{m}$.
Applying Theorem \ref{thm:av}, the classical Brunn-Minkowski inequality, and Lemma \ref{lem:pBM_balls},  we get
\begin{align*}
\vol\left(\frac{\partial K_{1}+\dots+\partial K_{m}}{m}\right)^{\frac{1}{n}} & =\frac{2}{m}\vol\left(\frac{\left(K_{1}+\dots+K_{m-1}\right)+\partial K_{m}}{2}\right)^{1/n}\\
 & \ge\frac{2}{m}\sqrt{\vol\left(K_{1}+\dots+K_{m-1}\right)^{1/n}\vol\left(K_{m}\right)^{1/n}}\\
 & \ge\frac{2}{m}\sqrt{\left(x_{1}+\dots+x_{m-1}\right)\,x_{m}}\\
 & \ge x_{1}^{\frac{1}{m}}\dots x_{m}^{\frac{1}{m}},
\end{align*}
that is 
\[
\vol\left(\frac{\partial K_{1}+\dots+\partial K_{m}}{m}\right)\ge x_{1}^{\frac{n}{m}}\dots x_{m}^{\frac{n}{m}}=\vol\left(K_{1}\right)^{\frac{1}{m}}\dots\vol\left(K_{m}\right)^{\frac{1}{m}},
\]
as claimed. 
Note that the inequality in Lemma \ref{lem:pBM_balls} is strict for $m>1$ and hence equality cannot occur in this case. 
\end{proof}

The proof of Theorem \ref{thm:BBM} is similar to that of Theorem \ref{thm:av}, but slightly more technical. We make use of the following fact. 
\begin{lemma}\label{lem:BBM_reduction}
	Let  $n\ge2$ and $\lambda\in[0,1]$. Let $R:(0,\infty)\to(0,\infty)$ be the function defined by
	$$
	R_n(x)=\frac{\left(|(1-\lambda)x+\lambda|^n-|(1-\lambda)x-\lambda|^n\right)\left(|(1-\lambda)+\lambda x|^n-|(1-\lambda)-\lambda x|^n\right)}{x^n}.
	$$
If $n=2$ then $R_2$ is a constant function.  If $n>2$ then  $R_n(x)>R(1)$ for all $x\in(0,\infty)$.
\end{lemma}

\begin{proof} Suppose first that $n=2$. Then 
\[|(1-\lambda)x+\lambda|^2-|(1-\lambda)x-\lambda|^2=|(1-\lambda)+\lambda x|^2-|(1-\lambda)-\lambda x|^2=4\lambda(1-\lambda)x \]
which means that $R_n(x)= 16\lambda^2(1-\lambda)^2 = R(1)$ is a constant function.  

Suppose next that $n>2$. Since the function $R_n$ remains the same when replacing $\lambda$ by $1-\lambda$,  we may assume without loss of generality that $0\le \lambda\le1/2$. For $\lambda=0$ the statement is trivial, so we  assume that  $0<\lambda\le1/2$.
		
		Denote the restrictions of $R_n$ to the intervals $(0,\frac{\lambda}{1-\lambda})$, $(\frac{\lambda}{1-\lambda},\frac{1-\lambda}{\lambda})$ and $(\frac{1-\lambda}{\lambda},\infty)$ by $f,g$ and $h$, respectively.  We shall show that $g$ attains a strict minimum at $x=1$, $h$ is non-decreasing and $f$ is non-increasing, from which it follows that $R_n(x)$ attains a strict  minimum at $x=1$.
		
		Let us first show that $g(x)>g(1)$ for all $x\in(\frac{\lambda}{1-\lambda},\frac{1-\lambda}{\lambda})$. Indeed, we have 
		\begin{equation}\label{eq:g}
		g(x)=\frac{\left(((1-\lambda)x+\lambda)^n-((1-\lambda)x-\lambda)^n\right)\left(((1-\lambda)+\lambda x)^n-((1-\lambda) -\lambda x)^n\right)}{x^n}
		\end{equation}
		which is well-defined for all $x\in(0,\infty)$ and identifies with $R_n(x)$ on  $x\in(\frac{\lambda}{1-\lambda},\frac{1-\lambda}{\lambda})$.  We   show that $\min_{x\in(0,\infty)} g(x)=g(1)$. To this end, note that $g$ is of the form  $$g(x)=a_{-n}x^{-n}+\dots+a_nx^n$$ for some non-negative reals $a_{-n},\dots,a_n$, which are not all zeroes. In fact, it is not hard to check that $g(x)$ is not a constant function and so there exists some $k\neq 0$ for which $a_k=a_{-k}>0$. Indeed, this claim immediately follows from the fact that the numerator in \eqref{eq:g} is a polynomial of degree at least $2(n-1)$ whereas the denominator in \eqref{eq:g} is the monomial $x^n$ of degree $n$, and the fact that $2(n-1)>n$ as $n>2$.

  For each $k\in\{-n,\dots,n\}$ denote $\lambda_k=\frac{a_k}{\sum_{k=-n}^{n}a_k}$.  Observing that $g(x)=g(1/x)$ for all $x\in(0,\infty)$, it follows that  $a_{-k}=a_{k}$ and hence $\lambda_k=\lambda_{-k}$ for all $k\in\{1,\dots,n\}$. By the arithmetic-geometric means inequality, we obtain that for all $x\in(0,\infty)$,
		$$
		g(x)=(\sum_{k=-n}^{n}a_k)\sum_{i=-n}^n\lambda_k x^k\ge (\sum_{k=-n}^{n}a_k)\prod_{k=-n}^{n}(x^k)^{\lambda_k}=\sum_{k=-n}^{n}a_k=g(1),
		$$
		Note that if $x\neq 1$, the above inequality is strict as $x^k\neq x^{-k}$ for some $k\neq0$ such that $a_k=a_{-k}>0$. Therefore, $g(x)>g(1)$ for all $x\in (\frac{\lambda}{1-\lambda}, \frac{1-\lambda}{\lambda})$, as claimed.
		
		Next, let us show that $h:(\frac{1-\lambda}{\lambda},\infty)\to(0,\infty)$ is non-decreasing. Indeed, we have 
		$$
		h(x)=\frac{\left(((1-\lambda)x+\lambda)^n-((1-\lambda)x-\lambda )^n\right)\left((\lambda x+(1-\lambda))^n-(\lambda x -(1-\lambda))^n\right)}{x^n}.
		$$
		For $\gamma\in(0,1)$ define the functions  $a_\gamma(x)=\gamma x + (1-\gamma)$ and $b_\gamma(x)=\gamma x - (1-\gamma)$. An easy calculation yields that 
		$$
		h'(x)=nx^{-(n+1)}c(x)\cdot\left(-1+\lambda x\frac{a^{n-1}_{\lambda}(x)-b^{n-1}_{\lambda}(x)}{a^n_{\lambda}(x)-b^n_{\lambda}(x)}+(1-\lambda)x\frac{a^{n-1}_{1-\lambda}(x)-b^{n-1}_{1-\lambda}(x)}{a^n_{1-\lambda}(x)-b^n_{1-\lambda}(x)}\right)
		$$
		where $c(x)=\big(a^n_{\lambda}(x)-b^n_{\lambda}(x)\big)\big(a^n_{1-\lambda}(x)-b^n_{1-\lambda}(x)\big)\ge0$.  
	Noting that $(a_\gamma(x)+b_\gamma(x))/2=\gamma x$ for all $\gamma\in(0,1)$, it follows that 
		\begin{align*}
			\gamma x\cdot\frac{a_{\gamma}^{n-1}\left(x\right)-b_{\gamma}^{n-1}\left(x\right)}{a_{\gamma}^{n}\left(x\right)-b_{\gamma}^{n}\left(x\right)} & =\frac{a_{\gamma}\left(x\right)+b_{\gamma}\left(x\right)}{2}\cdot\frac{a_{\gamma}^{n-1}\left(x\right)-b_{\gamma}^{n-1}\left(x\right)}{a_{\gamma}^{n}\left(x\right)-b_{\gamma}^{n}\left(x\right)}\\
			& =\frac{1}{2}\left(1+\frac{b_{\gamma}\left(x\right)a_{\gamma}^{n-1}\left(x\right)-a_{\gamma}\left(x\right)b_{\gamma}^{n-1}\left(x\right)}{a_{\gamma}^{n}\left(x\right)-b_{\gamma}^{n}\left(x\right)}\right)\\
			& =\frac{1}{2}+\frac{1}{2}a_{\gamma}\left(x\right)b_{\gamma}\left(x\right)\frac{a_{\gamma}^{n-2}\left(x\right)-b_{\gamma}^{n-2}\left(x\right)}{a_{\gamma}^{n}\left(x\right)-b_{\gamma}^{n}\left(x\right)}\ge\frac{1}{2}
		\end{align*}
	Therefore, 
		$$
		 -1+ \lambda x\frac{a^{n-1}_{\lambda}(x)-b^{n-1}_{\lambda}(x)}{a^n_{\lambda}(x)-b^n_{\lambda}(x)}+(1-\lambda)x\frac{a^{n-1}_{1-\lambda}(x)-b^{n-1}_{1-\lambda}(x)}{a^n_{1-\lambda}(x)-b^n_{1-\lambda}(x)}\ge 0,
   $$
		and so $h'(x)\ge0$. 
  
  Note that $f(x)=h(1/x)$ for all $x\in(0,\frac{\lambda}{1-\lambda})$, and hence $f'(x)=-h'(1/x)/x^2\le0$. This completes our proof. 
\end{proof}

\begin{proof}[Proof of Theorem \ref{thm:BBM}]
We begin with the proof of the asserted inequality. By homogeneity, we assume without loss of generality that $1=\vol(T)\le\vol(K)$.  Denote  $x=\vol^{\frac{1}{n}}(K)$,
		$L_1=(\lambda K\ominus (1-\lambda) T)\cup((1-\lambda) T\ominus\lambda K)$ and  $L_2=(\lambda T\ominus (1-\lambda) K)\cup((1-\lambda) K\ominus \lambda T)$.
		
    For $L_1 \ne \emptyset$, either $\lambda K\ominus (1-\lambda) T=\emptyset$ and $ L_1-\lambda K\sub (1-\lambda) T$  or $(1-\lambda) T\ominus \lambda K=\emptyset$ and $L_1 -\lambda T\sub (1-\lambda) K$. It follows by Brunn-Minkowski  inequality \eqref{eq:BMadditive} that 
		\begin{equation*}
			\vol(L_1)\le|\lambda\vol(K)^{\frac{1}{n}}-(1-\lambda)\vol(T)^{\frac{1}{n}}|^n=|\lambda x-(1-\lambda)|^n.
		\end{equation*}
		(This statement is trivially true for $L_1= \emptyset$.)
		
		Applying Proposition \ref{prop:decomp-of-sum} and the Brunn-Minkowski inequality \eqref{eq:BMadditive} once more, we get
		\begin{equation}\label{eq:f1}
			\begin{aligned}
				\vol(\lambda\partial K+(1-\lambda)\partial T)&=\vol(\lambda K+(1-\lambda)T)-\vol(L_1)\\
				&\ge (\lambda\vol(K)^{\frac{1}{n}}+(1-\lambda)\vol(T)^{\frac{1}{n}})^n-\vol(L_1)\\
				&=|\lambda x+(1-\lambda)|^n-\vol(L_1)\\
				&\ge |\lambda x+(1-\lambda)|^n- |\lambda x-(1-\lambda)|^n
			\end{aligned}
		\end{equation}
		Similarly, we get
		\begin{equation}\label{eq:f2}
			\vol(\lambda\partial T+(1-\lambda)\partial K)\ge |\lambda +(1-\lambda)x|^n -  |\lambda -(1-\lambda)x|^n
		\end{equation}
  As in Lemma \ref{lem:BBM_reduction}, denote
  \[
  R_n(x)=\frac{( |\lambda x+(1-\lambda)|^n- |\lambda x-(1-\lambda)|^n)(|\lambda +(1-\lambda)x|^n -  |\lambda -(1-\lambda)x|^n)}{x^n}.\]
		Combining  \eqref{eq:f1} and \eqref{eq:f2}, using the fact that  $\vol(K)\vol(T)=x^n$, and applying Lemma \ref{lem:BBM_reduction}, we have
		\begin{align*}
		\frac{\vol(\lambda\partial K+(1-\lambda)\partial T)\vol(\lambda\partial T+(1-\lambda)\partial K)}{\vol(K)\vol(T)}\ge R_n(x)\ge R_n(1)
		\end{align*}
		which completes the proof of the asserted inequality.

Next, we analyze the equality cases when  $K$ and $T$ are convex bodies and  $\lambda\neq 1/2$. 
%
%
First note that the equality case in the first inequality of  \eqref{eq:f1}, namely      \[\vol(\lambda K + (1-\lambda)T)= \lambda\vol(K)^\frac 1n + (1-\lambda)\vol(T)^\frac{1}{n},\]
means that the convex bodies $K,T$ are homothetic. 

Second, note that since $K$ and $T$ are homothetic, $\vol(K)\ge\vol(T)$ and $\lambda\neq 1/2$, it follows that  either $L_1\neq\emptyset$ or $L_2\neq\emptyset$. By symmetry in $\lambda$, we assume without loss of generality that $L_1\neq\emptyset$. Therefore, the equality case in the second inequality of \eqref{eq:f1}, namely
%
    \[\vol(L_1)=|\lambda\vol(K)^{\frac{1}{n}}-(1-\lambda)\vol(T)^{\frac{1}{n}}|^n\]
 must hold, which means that  $K$, $-T$ and $L_1$ are homothetic (up to translation and taking the closure of $L_1$). In particular, $K$ (and $T$) must be centrally-symmetric up to translation. 

It is left to analyze the equality cases for which $R_n(x)=R_n(1)$. Clearly, if  $x=1$ then $\vol(K)=\vol(T)$ which, in turn, means that $T$ is a translation of $K$. If $x\neq 1$ then by  Lemma \ref{lem:BBM_reduction} it must hold that $n=2$ and $R_2$ is the constant function $R_2(1)$.  This completes our proof. 
\end{proof}
\begin{remark}
   We mention that the definition of Minkowski subtraction as well as parts of our propositions are  connected to an area in robotics   
   where $K$ is the so called ``table'' (which need not be convex, and can contain ``obstacles'') on which a robot is allowed to move, and $-T$ signifies the robot itself which has a certain shape, and thus $K\ominus T$ are the ``allowed'' points on which one may place the robot without hitting the obstacles. 
   The paper \cite{Baram_2015} deals with questions of this form, where the notion of a ``hole'' $H$ having a trace or not is discussed. In particular they obtain a result similar to our \eqref{eq:decomposition2} stating that when computing the Minkowski sum of two sets with ``holes'', at least one of the sets can be ``completed'' without changing the sum. 
   In our setting most of these issues reduce   to the question of whether $\partial H-\partial Q$ has any ``holes'' in it, namely (considering Proposition \ref{prop:decomp-of-sum}), whether $H\ominus (-Q)$ is empty or not.
\end{remark}

\section{Restricted Minkowski Sums}\label{sec:restrict}

Since our main results in this paper concern the Minkowski addition of {\em subsets} of two given convex sets, they fit within the framework of ``restricted Minkowski sums'' suggested by Szarek and Voiculescu \cite{SV}  (see also \cite{AGM2}*{Section 4.5.3}). To this end, given two sets $A,B\sub\RR^n$ and a subset $\Theta \sub A\times B$, one defines their $\Theta$-restricted Minkowski sum
by 
\[ A +_\Theta B = \{ x+y: (x,y)\in \Theta\}.
\]
Szarek and Voiculescu \cite{SV} proved   that if $\vol(\Theta)$ is sufficiently close to $\vol(A\times B)$ then  an Arithmetic  Brunn-Minkowski inequality with exponent $2/n$ holds.

\begin{theorem}[S.~J.~Szarek, D.~Voiculescu]\label{thm:restriceted sums}
	Let $\rho \in (0,1)$, $n\in \mathbb{N}$ and let $A,B\sub \R^n$ be such that 
	\[\rho \le \left(\frac{\vol(B)}{\vol(A)}\right)^{\frac 1n}\le \rho^{-1}.\]
	Let $\Theta\sub A\times B\sub \R^{2n}$ be such that
	\[\vol(\Theta)\ge (1-c\min\{\rho\sqrt{n},1\})\vol(A)\vol(B).\]
	Then
	\[\vol(A+_\Theta B)^{2/n}\ge \vol(A)^{2/n}+\vol(B)^{2/n}.\]
	Here $c>0$ is a universal constant, independent of $\rho,n,A$ and $B$.
\end{theorem}

It turns out that Theorem \ref{thm:restriceted sums} can be applied, together with 
Proposition 
\ref{prop:decomp-of-sum}, to prove an arithmetic version for our main theorem, with   exponent $2/n$, under an extra volume ratio  condition on  the two bodies. 
\begin{theorem}\label{thm:restricted_BM}
	There exists  $n_0\in {\mathbb N}$ such that for all $n\ge n_0$ and  any $K,T\sub \R^n$ satisfying 
	\[\frac{1}{\sqrt{n}}\le \left(\frac{\vol( K)}{\vol (T)}\right)^\frac 1n \le \sqrt{n},\]
	it holds that
	\[\vol(\partial K+\partial T)^{2/n}\ge \vol(K)^{2/n}+\vol(T)^{2/n}.\]
	In fact, the theorem holds for all $n\ge 2$ if we replace the volume ratio condition by 
	$ {C}/{\sqrt{n}}\le \left({\vol( K)}/{\vol (T)}\right)^{1/n} \le \sqrt{n}/C$ for some universal $C>0$. 
\end{theorem}

\begin{remark}
	The volume ratio condition cannot be omitted as the following simple example shows: consider he pair $K$ and $aK$ with $a\to 0$. The left hand side is the clearly $O(a)$ as $a\to 0$  whereas the right hand side is equal to $(1+a^2)\vol(K)^{1/n}$ and in particular does not converge to $0$ as $a\to 0$. 
\end{remark}

\begin{proof}[Proof of Theorem \ref{thm:restricted_BM}]
	Assume without loss of generality that $\vol(K)\ge \vol(T)$, so that 
	\[ 	K + T  = (\partial K + \partial T ) \cup (K\ominus T).\]  Define $\Theta \subseteq K\times T\sub \R^{2n}$ as follows
	\[\Theta= \{(x,y) \in K \times T: x \not\in (K\ominus T-y)\}.\]
	Clearly this means that 
	\begin{equation}\label{eq:set comprison}
		K+_\Theta T\subseteq \partial K+\partial T.
	\end{equation}
	To apply Theorem \ref{thm:restriceted sums}  we need to find an appropriate lower bound for $\vol(\Theta)$.  First,  for any $y\in T$ we have  \[ \vol(K\setminus(K\ominus T-y))\ge \vol(K)-\vol(K\ominus T-y) = 
	\vol(K)-\vol(K\ominus T),
	\] 
	which implies the inequality  
	\begin{equation}\label{eq:set minus}
		\vol(\Theta)\ge \vol(T)(\vol(K)-\vol(K\ominus T)).
	\end{equation}
	Next, since $(K\ominus T) - T \sub K$, the classical Brunn-Minkowski inequality implies that 
	\begin{equation}\label{eq:minus bound}
		\vol(K\ominus T)^{1/n} \le \vol(K)^{1/n}-\vol(T)^{1/n}\le (1-r_n)\vol (K)^{1/n}
	\end{equation}
	where we have used the assumption in the statement of the proposition in the form $ (\frac{\vol( K)}{\vol (T)} )^\frac 1n \le r_n^{-1}$. 
	Combining the estimates \eqref{eq:set minus} and \eqref{eq:minus bound} we get that $\vol(\Theta) \ge (1-(1-r_n)^n)\vol(K)\vol(T)$, using the conditional bound on $\Theta$ from Theorem \ref{thm:restriceted sums} we can see that it is sufficent to find $r_n$ such that
	\[(1-(1-r_n)^n)\vol(K)\vol(T) \ge (1-c\min \{r_n \sqrt n,1\})\vol (K)\vol (T),\]
	where $c$ is the constant from Theorem \ref{thm:restriceted sums}. 
	Assuming (as we may) that $r_n\ge \frac{1}{\sqrt n}$ the inequality can be simplified to
	\[ r_n \ge 1-c^{1/n},\]
	meaning that by choice of $r_n\ge \max\{\frac{1}{\sqrt n}, 1-c^{1/n}\}$ we get the desired result. For $n\ge c_0$ this is simply $r_n \ge 1/\sqrt{n}$. 
\end{proof}
\bibliographystyle{amsplain}
\bibliography{ref}

\end{document}